%% file: bhf.tex
\documentclass[a4paper,12pt]{article}%{pcms2-l}
\usepackage{theorem,amsmath,latexsym,amssymb,newcom,esint,enumerate}
\usepackage[mathscr]{eucal}
\usepackage[dvips]{graphicx,color}
\begin{document}
\title{\textsf{A boundary partial regularity \\
               and a regularity criterion for \\
               Harmonic Heat flows}}
\author {\textsf{Kazuhiro HORIHATA}}
\maketitle
\begin{abstract}
{\bf{abstract}.} \enspace 
In my previous paper I have contrived a Ginzburg-Landau heat flow with a time-dependent parameter
and by using it,
I constructed a harmonic heat flow into spheres 
with a monotonical inequality and a reverse Poincar\'{e} inequality.
This paper establishes these two energy inequalities near the boundary
and then by making the best of them, we discuss a partial boundary regularity.
In addition to it, we demonstrate a whole domain's regularity
under \lq\lq{the one-sided condition.}\rq\rq 
This has been proposed by S.Hildebrandt and K.-O.Widman. %~\cite{hildebrandt-widman}.


\end{abstract}
\input intro_bhf.tex

\input glhf_bhf.tex
\input estimate_bhf.tex
\input bhf_bhf.tex
\input main_bhf.tex
\input ref_bhf.tex

\end{document}

%% file: intro_bhf.tex
%
%#! platex bhf
%
\setcounter{chapternumber}{1}\setcounter{equation}{0}
\renewcommand{\theequation}%
           {\thechapternumber.\arabic{equation}}
\section{\enspace Introduction.}
\label{SEC:Intro}
A previous paper: K.Horihata~\cite{horihata} constructed 
a harmonic heat flow between $d$-dimensional unit ball and $D$-dimensional unit sphere with
\renewcommand{\theenumi}{\roman{enumi}}\begin{enumerate}
\item
a global energy inequality
\item
a monotonical inequality
\item
a reverse Poincar\'e inequality.
\end{enumerate}
Furthermore it showed that the flow is smooth except a small singular set.
This paper first discusses a partial regularity near the boundary for such a flow 
by showing (ii) and (iii) and next establish the whole domain's smoothness
under a suitable condition for an initial maps.
\par 
We rigorously state our problem:
Let $d$ and $D$ be positive integers greater than or equal to $2$ and $T$ a positive number.
Suppose that $\Omega$ is a convex domain with $C^2$-boundary.
Set $Q(T)$ $=$ $(0,T)\times\Omega$. 
Give a mapping $u_0 \, \in \,$ $H^{1,2} (\Omega ; \mathbb{S}^D)$. We then consider the heat flow:
\begin{equation}
\left\{
\begin{array}{rll}
%1
\dfrac {\partial u}{\partial t} & \, = \,\triangle u \, + \, | \nabla u |^2 u
\quad & \mathrm{in} \;\; Q(T),
\\[5pt]
%2
u(0,x) & \, = \, u_0(x) \qquad & \mathrm{at} \quad \{0\} \times \Omega,
\\[3pt]
%3
u(t,x) & \, = \, u_0(x) \qquad & \mathrm{on} \quad [0,T) \times \partial \Omega.
\end{array}
\right.
\label{EQ:HHF}
\end{equation}
Define a function class:
\allowdisplaybreaks\begin{align}
%1
%& L^\infty \bigl(0,T;H^{1,2}(\Omega; \mathbb{S}^D)\bigr) 
%2
& V(Q(T);\mathbb{S}^D)\;:=\;L^\infty \bigl(0,T;H^{1,2}(\Omega;\mathbb{S}^D)\bigr)\cap
H^{1,2}\bigl(0,T;L^2(\Omega;\mathbb{R}^{D+1})\bigr).\notag\end{align}
The weak formulation of \eqref{EQ:HHF} is as follows:
For any given mapping $u_0 \, \in \,$ $H^{1,2} (\Omega;\mathbb{S}^D)$,
we call a mapping $u \, \in \, V (Q(T);\mathbb{S}^D)$
a weakly harmonic heat flow {(\textit{WHHF})} provided that
\allowdisplaybreaks\begin{alignat}{2}
%1
&
\lint_{Q(T)}\Bigl(\la \frac {\partial u}{\partial t},\phi\ra
\,+\,\la\nabla u, \nabla\phi\ra\,-\,\Bigr. && \Bigl. 
\la u,\phi\ra |\nabla u|^2\Bigr)\, dz\,=\,0,\label{EQ:1}\\
%2
& 
u(t)\,-\,u_0\in\tc(\Omega ;\mathbb{R}^{D+1})&&\quad \text{for almost every} \; t \, \in \, (0,T),
\label{EQ:2}\\
%3
&\underset {t \searrow 0}{\lim} \, u(t) \, = \, u_0
\quad && \quad \mathrm{in} \quad L^2 (\Omega ;\mathbb{R}^{D+1})
\label{EQ:3}\end{alignat} holds for any $\phi \in C_0^\infty ({Q(T)};\mathbb{R}^{D+1})$.
\par
K.Horihata~\cite{horihata} has proposed a new approximate evolutional scheme
said to be the Ginzburg-Landau heat flow (\textit{GLHF}):
\begin{equation} \dfrac{\partial u_\lambda}{\partial t}\,-\, \triangle u_\lambda
\,+\,\lambda^{1-\kappa}\bigl(|u_\lambda |^2\,-\,1\bigr) u_\lambda\;=\;0\quad\mathrm{in}\quad Q(T).
\label{EQ:GLHF}\end{equation}
Here notice that an application of a maximal principle on an unbounded domain 
referred to F.John~\cite[Chap 7]{john}
prevails $| u_\lambda |$ $\le$ $1$ in $Q(T)$.
The benefits of my scheme is to easily derive
\begin{equation}
\lint_{Q(T)} \lambda^{1-\kappa}( | u_\lambda |^2 \, - \, 1 )^2 \, dz
\;= \; O(1/\log\lambda)\quad\mathrm{as}\;\lambda\nearrow\infty.
\label{INEQ:GLQF0}\end{equation}
\par
We outlook history on a boundary regularity of harmonic mappings or harmonic heat flows:
For the energy minimizing map, R.Schoen and K.Uhlenbeck~\cite{schoen-uhlenbeck} has studied
a boundary regularity under a suitable condition.
A similar result can be found in J.Jost and M.Meier~\cite{jost-meier}.
J.Qing~\cite{qing} proved the boundary regularity of weakly harmonic maps 
from surfaces with the boundary.
C.Poon~\cite{poon} constructed a smooth harmonic map between $B^3$ and $\mathbb{S}^2$ 
except a prescribed point $a$ $\in$ $\overline{B}^3$,
noncontinuous at $a$ and $u_a$ $=$ $x$ on $\partial B^3$.
\par
On the other hand, there are few papers 
that discussed boundary regularity results on harmonic heat flows.
We refer it to Y.Chen \cite{chen-91} or Y.Chen-F.-H.Lin \cite{chen-lin} 
or Y.Chen-J.Li-F.H.Lin \cite{chen-li-lin} or C.Y.Wang \cite{wang-99}.
Since the map by C.Poon \cite{poon} is also a non-smooth weakly harmonic heat flow,
it may be reasonable to investigate a partial regularity for our harmonic heat flow.
\par
This paper has two folds: The first constructs a few energy inequalities on $\overline{\Omega}$ 
and then discusses a partial regularity for our harmonic heat flows near the boundary;
The proof identically proceeds to 
that of the author's previous paper of K.Horihata~\cite{horihata}
combined with that by Y.Chen~\cite{chen-91}.
The main technical ingredients are a monotonical inequality 
and reverse Poincar\'e inequality.
We present a simplified proof of the latter inequality.

After that by utilizing the monotonical inequality and the reverse Poincar\'e inequality,
we prove that the WHHF is smooth except on a small set.
More precisely, we assert
\begin{Thm}{\rm{(Partial Regularity Theorem).}}\label{THM:Main-1}
Let $d$ be a positive integer larger than $2$.
For a mapping $u_0$ $\in$ \begin{math}H^{1,2}(\Omega\,;\mathbb{S}^D)\end{math}
$\cap$ \begin{math}C^2(\Omega\setminus\Omega_{\delta_0})\end{math}
for a positive $\delta_0$ sufficiently small,
there exists a WHHF and it is smooth on a certain relative open set 
in $\overline{Q(T)}$ off a set % called {\bf{sing}}.
which has the finite $(d-\gamma_0)$-dimensional Hausdorff measure 
with respect to the parabolic metric, 
where $\gamma_0$ is a small positive number depending only on $u_0$ and $Q(T)$.
The WHHF also implies that for any \begin{math}\epsilon_0\,\in\,(0,1/2),\end{math} 
there exists a positive constant $C(\epsilon_0)$ such that
\allowdisplaybreaks\begin{align}
%1
& \lint_{t_0-(2R_1)^2}^{t_0-R_1^2}\lint_\Omega |\nabla u|^2 G_{z_0} \, dz
\,+\, 2\lint_{R_1}^{R_2}\, dR\lint_{t_0-(2R)^2}^{t_0-R^2}\, dt
\lint_\Omega
\Bigl\vert\frac{\partial u_\lambda}{\partial t}\,-\,\frac{x-x_0}{2\sqrt{t_0-t}}\cdot\nabla u\Bigr\vert^2\, G_{z_0}\, dx
\notag\\
%3
&
\;\le\;C \lint_{t_0-(2R_2)^2}^{t_0-R_2^2}\,dt\lint_\Omega|\nabla u |^2 G_{z_0} \, dx
\label{INEQ:Boundary-Monotonical-1}
\,+\,C(\epsilon_0) (R_2^{\epsilon_0}\,-\,R_1^{\epsilon_0})
\end{align}
holds for an arbitrary point $z_0$ $=$ $(t_0,x_0)$ in $\overline{Q}(T)$ and
for positive numbers $R_1$ and $R_2$ with $R_1 \le R_2$ satisfying $R_2$ $<$ $\sqrt{t_0/4}$,
where a function $G_{z_0}$ means a backward heat kernel indicated by
\begin{equation}
G_{z_0} (t,x)\;=\;\frac 1{\sqrt{4\pi (t_0-t)}^d}\exp\Bigl(-\frac{|x-x_0|^2}{4(t_0-t)}\Bigr).
\end{equation}
In addition to it, we obtain
\allowdisplaybreaks\begin{align}
%1
&\lint_{P_R (z_0)\cap Q(T)}|\nabla u|^2 \, dz
\;\le\;\frac{C}{R^2}\lint_{P_{2R}(z_0)\cap Q(T)}\vert u\,-\,a(t)\vert^2 \, dz\notag\\
%2
& \,+\, C\lint_{P_{2R}(z_0)\cap Q(T)}\vert \nabla u_0\vert^2 \,dz
\end{align}
for any parabolic cylinder $P_{2R}$ $(z_0)$ and
any \begin{math}a\,\in\,L^2(\mathbb{R}_+;\mathbb{R}^{D+1})\end{math},
where $C$ is a positive constant depending only on $Q(T)$ and $u_0$. 
\end{Thm} 
\par
The second establishes the smoothness of our heat flow under a certain condition
referred to M.Giaquinta~\cite[p.214, Theorem 1.1]{giaquinta} 
and M.Giaquinta and E.Giusiti~\cite[p.43, Theorem 5.2]{giaquinta-giusti}.
\begin{Thm}{\rm{(Regularity Criterion).}}\label{THM:Main-3}
Let $\Omega$ be compact with $C^2$-boundary and 
$u_0$ be a mapping of 
\begin{math}H^{1,2}(\Omega;\mathbb{S}^D)\end{math}
and \begin{math}C^2(\Omega\setminus\Omega_{\delta_0};\mathbb{R}^{D+1})\end{math}
for a positive number sufficiently small.
If the range of $u_0$, $u_0$ $(\Omega)$ is compactly contained in 
\begin{math}\mathbb{S}^D\cap\{y_d>0\}\end{math}
after a suitable rotation 
then the range of $u({Q}(T))$ is done there and they are smooth in $Q(T)$.
\end{Thm}
\par
In a forthcoming paper, we extend the results here
and the ones of K.Horihata~\cite{horihata} to a heat flow 
between any compact Riemannian manifolds.
\par
We adopt notation from K.Horihata~\cite{horihata} and list an additional notation:
We also added the new symbols mostly related to the boundary.
\vskip 9pt
\begin{center}\underline{\textit{Notation}}\end{center}
\vskip 9pt
\renewcommand{\labelenumi}{(\roman{enumi})}
\begin{enumerate}
\item
${Q(T)}=(0,T) \times \Omega$.
$\; \partial{Q(T)}\,=\,$ $[0,T)\times\partial\Omega$ $\cup$ $\{0\} \times \Omega$.
\item
$\mathbb{S}^D\,=\,$ $\{u=(u^1,u^2,\cdots,u^{D+1})\,\in \,\mathbb{R}^{D+1};$ 
$|u| \, = \, \sqrt{\sum_{i=1}^{D+1} (u^i)^2} \, = \, 1\}$.
\item
Set a points $x$ $=$ $(x_1,x_2,\cdots,x_d)$, $x_0$ $=$ $(x_{0,1},x_{0,2},\ldots,x_{0,d})$ and
$z_0\,=\, (t_0, x_0)\in\mathbb{R}^{d+1}$;
Indicate for any $r>0$, 
\allowdisplaybreaks\begin{align*}
%1
B_r(x_0)&\,=\,\{x\in\mathbb{R}^d\,;\, |x\,-\,x_0|<r\},\\ 
%2
P_r(z_0)&\,=\,(t_0-r^2,t_0+r^2)\times B_r(x_0), \; Q_r (z_0)\,=\,(t_0-r^2,t_0)\times B_r(x_0).
\end{align*}
In $B_r(x_0)$,$P_r(z_0)$ and $Q_r(z_0)$, 
the points of $x_0$ or $z_0$ will be often abbreviated when no confusion may arise.
\item
Letter $C$ denotes a generic constant depending only on $d$, $D$ and $u_0$.
By the letter $C(B)$ it means that the constant depends only on a parameter $B$.
\item
Let $A$ be a set in $\mathbb{R}^d$.
\begin{equation}A_{\delta_0}\;=\;\{x\,\in\, A\,;\,\dist(x,\partial A)>\delta_0\}\notag\end{equation}
for any positive number $\delta_0$.
\item
A vector $\nu$ and $\tau$ respectively denotes 
the unit outer normal and the tangential field along $\partial \Omega$.
\end{enumerate}

%% file: glhf_bhf.tex
9%
%#! platex bhf
%
\setcounter{chapternumber}{2}\setcounter{equation}{0}
\renewcommand{\theequation}{\thechapternumber.\arabic{equation}}
\section{\enspace GLHF.}
\label{P:GLHF}
The chapter introduces two fundamental energy inequalities:
The first is a monotonical inequality, 
Moreover the second is a hybrid inequality.
The proof of the latter is more straightforward than the one in the author's previous paper,
K.Horihata~\cite[Theorem 2.6]{horihata}.
\subsection{\enspace Monotonical Inequality.}
We introduce the usual monotonical inequality by Y.Chen and F.H.Lin~\cite{chen-lin}. % and
\begin{Thm}{\rm{(Monotonical Inequality).}}\label{THM:Mon}
There exist constants $C$ and $C(\mu_0)$ with $C(\mu_0)$ $\nearrow$ $\infty$ as $\mu_0$ $\searrow$ $0$
such that the following holds for any point $z_0 \, = \, (t_0,x_0)$ $\, \in \, Q (T)$ and 
for positive numbers $R$, $R_1$ and $R_2$ with $0 < R_1 < R_2$
with $\max(R_2,R)$ $\, \le \, \sqrt{t_0}/2$,
\allowdisplaybreaks\begin{align}
%1
&\lint_{t_0 - 4R_1^2}^{t_0 - R_1^2} \, dt \lint_{\Omega}\mathbf{e}_\lambda \, G_{z_0} \, dx      
\, + \, \lint_{R_1}^{R_2} \, dt \lint_{\Omega}\left| \frac {\partial u_\lambda}{\partial t}
\;-\;\frac {x\,-\,x_0}{2(t \, - \, t_0)} \cdot \nabla u_\lambda \right|^2 G_{z_0} \, dx
\notag\\
%2
&\; \le \; C \exp{(R_2^{\mu_0} \, - \, R_1^{\mu_0})} \lint_{t_0 - 4R_2^2}^{t_0 - R_2^2} \, dt \lint_{\Omega}
\mathbf{e}_\lambda \, G_{z_0} \, dx
\notag\\
%3
&\, + \, C(\mu_0) (R_2\,-\,R_1)\label{INEQ:Mon}\end{align}
with 
\allowdisplaybreaks\begin{align}
%1
& G_{z_0} (t,x) \; = \; \frac 1{\sqrt{4\pi (t_0-t)}^d}
\exp \Bigl( - \frac {|x-x_0|^2}{4(t_0-t)}\Bigr)
\quad ((t,x) \, \in \, (0,t_0) \times \mathbb{R}^d).
\notag\end{align}\end{Thm}
\subsection{\enspace Hybrid type Inequality for GLHF.}
We demonstrate the reverse Poincar\'{e} inequality.
To this end we prepare the below: Let the mapping $h_0$ be the solution to
\begin{equation}
\left\{
\begin{array}{rl}
- \triangle h_0 & \; = \; 0 \quad \mathrm{in} \quad \Omega, \\
h_0 & \; = \; u_0 \quad \mathrm{on} \quad \partial\Omega.
\label{EQ:HF}\end{array}\right.\end{equation}
Then we claim
\begin{Thm}{\rm{(Hybrid Inequality).}}\label{THM:HI}
For any positive number $\epsilon_0$ and any point $z_0$ in $Q(T)$, 
there exists a positive constant
$C(\epsilon_0)$ satisfying $C(\epsilon_0)$$\nearrow$ $\infty$ as $\epsilon_0 \searrow 0$
such that the inequality
\allowdisplaybreaks\begin{align}
%1
& \lint_{P_{R} (z_0) \cap Q(T)} \mathbf{e}_\lambda \, dz
\,\le\;\epsilon_0\lint_{P_{2R} (z_0)\cap Q(T)}\mathbf{e}_\lambda\, dz
\,+\,\frac {C(\epsilon_0)}{R^2}\lint_{P_{2R} (z_0)\cap Q(T)}\vert u_\lambda\,-\,h_{x_0}\vert^2\, dz
\notag\\
%2
& \,+\,\lint_{P_{2R}(z_0)\cap Q(T)} 
\vert\nabla h_0 (x)\vert^2\, dz \,+\,o(1)\quad ( \lambda \nearrow \infty)
\label{INEQ:HI}\end{align}
holds for any parabolic cylinder $P_{2R} (z_0)$.\end{Thm}
\vskip 9pt
\noindent{\underbar{Proof of Theorem \ref{THM:HI}.}}
\vskip 9pt\rm\enspace\par
\par Preliminary, we rewrite \eqref{EQ:GLHF} to it by the flatten-out coordinate:
To be more precise, fix $x_0$ $\in$ $\partial\Omega$ 
and choose an open set $U(x_0)$ and a $C^2$-function $\phi_{x_0}$ on
\begin{math}U(x_0)\,\cap\,\mathbb{R}^{d-1}\end{math},
where for any point
$x$ $=$ $(x',x_d)$ $=$ $(x_1,x_2,\ldots,x_d)$ $\in$ $U(x_0)$,
the following holds
\begin{equation}
\left\{\begin{array}{lc}
x\,>\,\phi_{x_0}(x')\;&\quad (x\,\in\,U(x_0)\cap\Omega)\\
x\,=\,\phi_{x_0}(x')\;&\quad (x\,\in\,U(x_0)\cap\partial\Omega)\\
x\,<\,\phi_{x_0}(x')\;&\quad (x\,\in\,U(x_0)\setminus\overline\Omega)
\end{array}\right.\end{equation}
and define a $C^2$-mapping $\Phi_{x_0}$ $=$ $(\Phi_{x_0}^1,\Phi_{x_0}^2,\cdots,\Phi_{x_0}^d)$
$:$ $U(x_0)$ $\to$ $\mathbb{R}^d$ by
\begin{equation}
\left\{\begin{array}{lc}
y_i\;=\;x_i\;=\;\Phi_{x_0}^i(x)\quad (i\,=\,1,2,\ldots,d)\\
y_d\;=\;x_d\,-\,\phi_{x_0}(x')\;=\;\Phi_{x_0}^d(x).
\end{array}\right.\end{equation}
Hereafter we use $y$ $=$ $(y',y_d)$ $=$ $(y_1,y_2,\ldots,y_d)$ $\in$ $\Phi_{x_0}(U(x_0))$.
\par
Choose $\overline{B_R} (0)$ $\subset$ $\Phi_{x_0} (U (x_0))$
and induce the mapping $v_\lambda$ and $h_{x_0}$ on $\overline{B_R}(0)$ by 
\begin{equation}
v_\lambda (t,y)\;=\; u_\lambda(t, \Phi_{x_0}^{-1}(y)),\;
h_{x_0} (y) \; = \; h_{0} (\Phi_{x_0}^{-1} (y)).
\end{equation}
and extend the mapping $v_\lambda \,- \, h_{x_0}$ in 
$\overline{B_R}(0)\cap\mathbb{R}_+^d$ to the one in $\overline{B_R}(0)$ 
designated by 
\begin{equation}\left\{\begin{array}{rlc}
%1
(v_\lambda\,-\,h_{x_0})(t,(y^\prime,y_d))&
\quad\mathrm{if}\quad y\,\in\,\overline{B_R} (0)\cap\mathbb{R}_+^d,\notag\\[2mm]
%2
0&
\quad\mathrm{if}\quad y\,\in\,\overline{B_R}(0) \cap \mathbb{R}_0^d,\notag\\[1mm]
%3
-( v_\lambda \, - \, h_{x_0} ) (t,(y^\prime,-y_d)) &
\quad\mathrm{if}\quad y\,\in\,\overline{B_R}(0)\cap\mathbb{R}_-^d
\notag\end{array}\right.\label{EQ:GLHF-Flat}\end{equation}
\allowdisplaybreaks\begin{align}
%1
& D_i \; = \; 
\left\{\begin{array}{lc}
\nabla_i\,-\,a_i\nabla_d& \quad(y_d\,>\,0)
\notag\\[2mm]
0& \quad(y_d\,=\,0)\notag\\[2mm]
\nabla_i\,+\,a_i\nabla_d&\quad (y_d\,<\,0)\\
\end{array}\right.
\;(i\,=\,1,2,\ldots,d-1)
\notag\\
%2
& \mathrm{with}\quad
a_i\;=\;\frac{\partial\phi_{x_0}}{\partial y_i}\quad (i\,=\,1,2,\ldots,d-1),\notag\\
%3
&
D_d\;=\;\nabla_d, \;
%4
D\;=\;(D_1,D_2,\cdots,D_d),\; L \; = \; \sum_{i=1}^{d-1} D_i^2 \, + \, D_d^2,\notag\\
%5
& D_\nu\;=\;\sum_{i=1}^d \frac {y^i}{|y|}D_i,\; D_\tau\;=\;D\,-\,\frac {y}{|y|}\,D_\nu ,\notag\\
%6
& \triangle_\tau\;=\;\triangle\,-\, 
\frac 1{\rho^{d-1}}\frac{\partial}{\partial \rho}
\Bigl(\rho^{d-1}\frac{\partial}{\partial\rho}\Bigr).\end{align}
\par We denote the extensive mapping by the same symbol $v_\lambda \, - \, h_{x_0}$:
Note that the mapping $v_\lambda \, - \, h_{x_0}$
belongs to $C^1$ $(B_R (0) \, ; \, \mathbb{R}^{D+1})$
and it satisfies the following identity:
\begin{eqnarray}
%1
\lint_{B_R(0)}\Bigl\langle \dfrac{\partial v_\lambda}{\partial t},\phi\Bigr\rangle \, dy
\, + \, \lint_{B_R(0)}\langle D (v_\lambda - h_{x_0} ), D \phi \rangle \, dy
\notag\\
%2
\, + \, \lambda^{1-\kappa} \lint_{B_R(0)}\bigl( | v_\lambda |^2 \, - \, 1 \bigr) 
\langle v_\lambda, \phi \rangle \, dy\; = \; 0
\label{EQ:GLHF-Boundary}\end{eqnarray}
for any mapping $\phi$ $\in$ $C_0^\infty$ $(B_R (0)\,;\,\mathbb{R}^{D+1})$.
Here we used
\allowdisplaybreaks\begin{align}
%1
&     
\langle D (v_\lambda - h_{x_0} ), D\phi\rangle 
\;=\;\sum_{i=1}^d \sum_{j=1}^{D+1} D_i (v_\lambda - h_{x_0} )^j D_i \phi^j
\notag\\
%2
&
\Bigl\langle \frac{\partial v_\lambda}{\partial t},
\phi\Bigr\rangle\;=\;\sum_{j=1}^{D+1}\frac{\partial v_\lambda^j}{\partial t}\phi^j,
\;\langle v_\lambda,\phi\rangle\;=\;\sum_{j=1}^{D+1} v_\lambda^j\phi^j. \notag\end{align}
\par Secondary, we take up a few inequalities:
\begin{Lem}\label{LEM:GLQE}
For any balls $B_{\rho_1}$ and $B_{\rho_2}$ 
with $B_{\rho_1}$ $\, \subset \,$ $B_{\rho_2}$ $\, \subset \,$ $B_R (0)$ and
any cylinders $P_{\rho_1}$ and $P_{\rho_2}$ with $P_{\rho_1}$ $\, \subset \,$ $P_{\rho_2}$
$\,\subset\,$ $P_R (0)$, the following holds{\rm{:}}
\allowdisplaybreaks\begin{align}
%1
& \lint_{P_{\rho_1}}\lambda^{1-\kappa}(1\,-\, |v_\lambda |^2 )\, dz
\;\le\;C\lint_{P_{\rho_2}}\mathbf{e}_\lambda\,dz\label{INEQ:GLQE}\\
%2
& \,+\,\frac {C}{(\rho_2-\rho_1)^2}\lint_{P_{\rho_2}} (1\,-\,|v_\lambda|^2)\, dz.
\notag\end{align}
\end{Lem}
\par Let $r$ be a positive number with \begin{math}R\,\le\,r\,<\,2R.\end{math}
Combined Theorem \ref{THM:Mon} with Lemma \ref{LEM:GLQE} above, we claim
\begin{Rem}\label{REM:DensityOfGLQE}
\begin{equation}
\sup_{0<r<2R} \frac 1{r^d} \lint_{P_r\cap{Q(T)}} \lambda^{1-\kappa}(1\,-\,|v_\lambda|^2)\, dz\;<\; C\,<\,\infty
\label{INEQ:DensityOfGLQE}\end{equation}
where $C$ is a positive constant independent of $R$.
\end{Rem}
Next we mention Poincar\'{e} type inequality:
\begin{Lem}\label{LEM:Poincare}
For any positive number $\epsilon_0$, 
there exists a positive constant $C(\epsilon_0)$ 
with $C(\epsilon_0)\,\nearrow\,\infty$ as $\epsilon_0\,\searrow\,0$
suct that the following holds\rm{:}
\allowdisplaybreaks\begin{align}
%1
& \sup_{t\,\in\,(t_0-r^2,t_0+r^2)}\fint\limits_{\partial B_r (x_0)} |v_\lambda\,-\,h_{x_0}|^2 \,d\mathcal{H}_y^{d-1}
\label{INEQ:Poincare}\\
%2
&
\;\le\;\frac{\epsilon_0^2}{r^d}\lint_{P_{2r}(z_0)} e_\lambda \, dz
\,+\, \frac 1 {r^{d-2}} \lint_{P_{2r}(z_0)}
\left\vert \frac{\partial v_\lambda}{\partial t}\right\vert^2 \, dz
\,+\,\frac{C(\epsilon_0)}{r^{d+2}}\lint_{P_{2r}(z_0)}\left\vert v_\lambda\,-\,h_{x_0}\right\vert^2 \, dz
\notag\end{align}
for any \begin{math}P_{2r}(z_0)\;=\;(t_0-(2r)^2,t_0+(2r)^2)\times B_{2r}(x_0)\end{math}
\begin{math}\;\subset\;\mathbb{R}^{d+1}.\end{math}\end{Lem}
\par Final step of the preparation lists symbols and auxiliary functions employed only here.
Give $L_\lambda$ by $[\log(\lambda^3/g(\lambda))]+1$
where a function $g(\lambda)$ $(\lambda\in\mathbb{R})$ is positive
and $g(\lambda)$ $\to$ $0$ $(\lambda\nearrow\infty)$.
Put $r$ $=$ $r(t)$ be any positive function satisfying the following:
\begin{Def}\label{DEF:Average}
For R given above, choose a function $r$ $=$ $r(t)$ with
\allowdisplaybreaks\begin{align}
%1
& r\lint_{\partial B_r^+}|\nabla (v_\lambda-h_{x_0})|^2 \, d\mathcal{H}_{y}^{d-1}
\,+\, r^3\lint_{-(2R)^2}^{-(2R)^2}\, dt\lint_{\partial B_r^+}
\left|\frac{\partial v_\lambda}{\partial t}\right|^2\, d\mathcal{H}_{y}^{d-1}\notag\\
%2
& \,+\, r\lint_{-(2R)^2}^{-(2R)^2}\lambda^{1-\kappa}\, dt\lint_{\partial B_r^+}
(|v_\lambda|^2-1)^2\, d\mathcal{H}_{y}^{d-1}\notag\\
%3
&\;\le\;2\lint_{B_{2R}^+\setminus B_R^+} |\nabla (v_\lambda-h_{x_0})|^2 \, dy
\,+\, 8R^2\lint_{-(2R)^2}^{-(2R)^2} \, dt \lint_{B_{2R}^+\setminus B_R^+}
\left|\frac{\partial v_\lambda}{\partial t}\right|^2 \, dy \notag\\
%4
& \,+\, \lint_{-(2R)^2}^{-(2R)^2}\lambda^{1-\kappa}\, dt\lint_{B_{2R}^+\setminus B_R^+}
(|v_\lambda|^2-1)^2\, dy.\label{INEQ:Average}\end{align}\end{Def}
\par We then introduce the decompositional convention:
\allowdisplaybreaks\begin{align}
%1
&\triangle\tau_0\;=\;1,\triangle\tau_l \;=\;(1/2)^l,\;\triangle \check{r}_l\;=\; Cr\triangle\tau_l,
\;\triangle\widetilde{r}_l\;=\;\widetilde{C} r\triangle\tau_l^2,
\notag\\
%2
& t_{l,m}\;=\;mr^2\triangle\tau_l^2 ,\; I_{l,m}\;=\;[t_{l,m-1},t_{l,m})
\quad (l\,=\,1,2,\ldots,L_\lambda\,;\,m\,=\,1,2,\ldots,1/\triangle\tau_l^2)\notag
\end{align}
with \begin{math}C^{-1}\,=\, \sum_{l=1}^{L_\lambda}\triangle\tau_l\end{math} and
\begin{math}\widetilde{C}^{-1}\,=\,\sum_{l=1}^{L_\lambda}\triangle\tau_l^2,\end{math} 
\allowdisplaybreaks\begin{align}
%1
{\Check{r}}_1\;=\;& (1\,-\,\epsilon_0^4)r,\;
{\Check{r}}_l\;=\;{\Check{r}}_1\,+\,\epsilon_0^4\sum_{j=1}^{l-1}\triangle r_j,\;
\notag\\
%2
{\widetilde{r}}_1\;=\;& (1\,-\,\epsilon_0^4)r,\;
{\widetilde{r}}_l\;=\;{\widetilde{r}}_1
\,+\,\epsilon_0^4\sum_{j=1}^{l-1}\triangle\widetilde{r}_j\;
(l\,=\,2,3,\ldots, L_{\lambda})\notag. \notag\end{align}
There exists a positive number $r_l$ with $r$ $<$ $r_l$ $<$ $2r$ with
\allowdisplaybreaks\begin{align}
%1
&  
\lint_{t_{l,m-1}}^{t_{l,m}}\, dt 
\lint_{B_{\check{r}_l+\triangle{r}_l/4}\setminus B_{\check{r}_l-\triangle{r}_{l+1}/4}}
|Dv_\lambda (t,y)|^2 \, dy
\;=\;\frac 4{3\triangle r_l}\lint_{t_{l,m-1}}^{t_{l,m}} \, dt
\lint_{\partial B_{r_l}} |Dv_\lambda (t,y)|^2 \,d\mathcal{H}_y^{d-1}.\notag\end{align}
Give
\allowdisplaybreaks\begin{align}
%1
& B_0\;=\;\emptyset, B_l\;=\;B_{r_l},\;\widetilde{B}_{l}\;=\;B_{\widetilde{r}_{l}},\;
{T}_l\;=\;{B}_l\setminus{B}_{l-1}, 
{\widetilde{T}}_l\;=\;{\widetilde{B}}_l\setminus{\widetilde{B}}_{l-1}\notag\\
%3
& D^+\;=\;D\cap\{y_d>0\}\;\mathrm{for}\;\mathrm{any}\;\mathrm{set}\; D\;\subset\;\mathbb{R}^d
\quad (l\,=\,1,2,\ldots,L_\lambda).
\notag\end{align}
\par
Construct support function: 
Let \begin{math}(\rho,\omega_{d-1};t)\,\in\,(0,r)\times{S^{d-1}}\times(0,r^2);\end{math}
Then set \begin{math}{f}_{\lambda}\;=\;{f}_\lambda(\rho,\omega_{d-1})\;=\;{f}_{\lambda}(\rho,\omega_{d-1};t)\end{math}
by the solution to
\allowdisplaybreaks\begin{align}
%1
%\left\{
%\begin{array}{lll}
%1.1
&\Bigl(\dfrac{\partial}{\partial{\rho}}\,+\,\triangle_\tau\Bigr){f}_{\lambda}\;=\;0
\quad\mathrm{in}\quad (0,r)\,\times\,\mathbb{S}^{d-1}, 
\notag\\[2mm]
&{f}_{\lambda}\;=\;\left\{
\begin{array}{lll}
&({v}_{\lambda}\,-\,h_{0})(t,(y',y_d))&\; (y_d>0)\\
&0&\; (y_d= 0)\\
-&({v}_{\lambda}\,-\,h_{0})(t,(y',-y_d))&\; (y_d<0)
\end{array}\right.
\quad\mathrm{on}\quad \{r\}\,\times\,\mathbb{S}^{d-1}. \label{EQ:Support}\end{align}
\begin{Rem}\label{REM:Sym}
$f_\lambda (\rho,(y',0)/\rho)$ $\,\equiv\,0$ in 
\begin{math}(\rho,y')\,\in\, (0,r)\times\mathbb{S}^{d-1}\cap\{y_d=0\}.\end{math}
Indeed if we give $\tilde{f}_\lambda$ by
\allowdisplaybreaks\begin{align}
%1
& \left(\frac{\partial}{\partial\rho}\,+\,\triangle_\tau\right)\tilde{f}_\lambda\;=\;0
\quad (0,r)\times\mathbb{S}^{d-1},\notag\\
%2
& \tilde{f}_\lambda(r,(y',y_d)/r)\,=\,
\left\{\begin{array}{lll} 
-&(v_\lambda\,-\,h_{x_0})(t,r,(y',y_d)/r)& (y_d> 0)\\
 &0& (y_d=0)\\
 &(v_\lambda\,-\,h_{x_0})(t,r,(y',-y_d)/r)& (y_d<0)
\end{array}\right.
\;\mathrm{on}\quad y\,\in\,\mathbb{S}^{d-1}.
\notag\end{align}
Thus we have
\begin{equation*}
{f}_\lambda\,+\,\tilde{f}_\lambda\;=\;0\;\mathrm{at}\;\{r\}\times\mathbb{S}^{d-1}
\;\mathrm{and}\;
\left(\frac{\partial}{\partial\rho}\,+\,\triangle_\tau\right)(f_\lambda+\tilde{f}_\lambda)\;=\;0
\;\mathrm{in}\; (0,r)\times\mathbb{S}^{d-1}
\end{equation*}
which implies \begin{math}{f}_\lambda\,+\,\tilde{f}_\lambda\;\equiv\;0.\end{math}
Since 
\begin{equation*}{f}_\lambda(\rho,(y',0)/\rho)\,=\,\tilde{f}_\lambda(\rho,(y',0)/\rho),\end{equation*}
we deduce
\begin{equation*}{f}_\lambda(\rho,(y',0)/\rho)\,=\,0.\end{equation*}\end{Rem}
\par
A first step establishing a Hybrid type inequality is to construct a certain support mapping
$w_{\lambda,l}$ by making the best of the function $f_\lambda$:
They are given by the solutions of
\allowdisplaybreaks\begin{align}
%1
&\left\{\begin{array}{rcl}
%1.1
& -\triangle{w}_{\lambda,l}\;=\;0 & \quad\mathrm{in}\quad T_l,\\
%1.2
& {w}_{\lambda,l}\;=\;f_{\lambda}(\widetilde{r}_l,\omega_{d-1})
& \quad\mathrm{on}\quad\partial B_l,\\
& {w}_{\lambda,l}\;=\;f_{\lambda}(\widetilde{r}_{l-1},\omega_{d-1})
& \quad\mathrm{on}\quad\partial B_{l-1}.
\end{array}\right.
\label{EQ:Support-4}\end{align}%\\
\par
As in the same way as in Remark \ref{REM:Sym}, we must remark
\begin{equation*}w_{\lambda,l} \;\equiv\;0\quad\mathrm{on}\; T_l^+\cap\{y_d\,=\,0\}\end{equation*} holds.
We state a few properties for the mapping $w_{\lambda,l}$ used below:
To introduce it, we shall recall a sequence of hyper spherical-harmonics $\{ \phi_n^{(\alpha)}\}$ 
$( n = 0, 1, \ldots ; \alpha = 1, 2, \ldots, N(n))$
where $\phi_n^{(1)}$, $\phi_n^{(2)}$, $\cdots,$ $\phi_n^{(N(n))}$ on $\mathbb{S}^{d-1}$
are independent components with degree $n$. 
\par
Writing \begin{math}y\,=\,(y^\prime,y_d)\,=\,(r,\omega_{d-1})\end{math}, 
we have the expression formula:
\allowdisplaybreaks\begin{align}
%1
&
f_{\lambda}^{n,(\alpha)}(r)
\;=\;\fint\limits_{{\mathbb{S}^{d-1}}}\langle f_{\lambda} (r,\omega_{d-1}),\phi_n^{(\alpha)}(\omega_{d-1})\rangle
\,d\omega_{d-1}\quad (r=\widetilde{r}_{l-1},\widetilde{r}_{l}),\notag\\
%2
& w_{\lambda,l} (x)\;=\; w_{\lambda,l} (r,\omega_{d-1})
\notag\\
%3
& \;=\;\sum_{n=1}^\infty\sum_{\alpha=1}^{N(n)} 
a_n^{(\alpha)} \Bigl( \frac r {{r_l}} \Bigr)^n \phi_n^{(\alpha)} (\omega_{d-1})
\,+\,\sum_{n=1}^\infty\sum_{\alpha=1}^{N(n)} b_n^{(\alpha)}
\Bigl( \frac {{r_{l-1}}} r \Bigr)^{n+d-2}\phi_n^{(\alpha)} (\omega_{d-1})
\notag\\
%4
& \,+\,a_0^{(1)}\label{EQ:Formula}\intertext{with}
%5
& a_n^{(\alpha)}\;=\;
\frac{f_{\lambda}^{n,(\alpha)} (\widetilde{r}_{l})
\,-\, 
{f}_{\lambda}^{n, (\alpha)} (\widetilde{r}_{l-1})\tau_l^{n+d-2}}
{1-\tau_l^{2n+d-2}}, 
%6
b_n^{(\alpha)}\;=\;-\frac {{f}_{\lambda}^{n,(\alpha)} (\widetilde{r}_l)\tau_l^{n}
\,-\,
{f}_{\lambda}^{n,(\alpha)}(\widetilde{r}_{l-1})}{1-\tau_l^{2n+d-2}}, \quad
\notag\\
%7
&\tau_l\;=\;\frac {{r_{l-1}}}{{r_l}}.
\notag\end{align}
A direct computation enjoys
\begin{Lem}\label{LEM:W-F}
$w_{\lambda,l}$ satisfies
\allowdisplaybreaks\begin{align}
%1
& \sum_{l=1}^{L_\lambda}\lint_{T_l}
(|\nabla_\nu w_{\lambda,l}|^2\,+\,|\nabla_\tau w_{\lambda,l}|^2) \, dy
%2
\;\le\;C r\lint_{\partial B_r}|\nabla_\tau (v_\lambda-h_{x_0})|^2 \, d\mathcal{H}_y^{d-1}.
\label{INEQ:F-W}\end{align}\end{Lem}

%% file: estimate_bhf.tex
%
%#! platex bhf
%
\vskip 2pt\par
After the preparation above, we demonstrate the proof of Theorem \ref{THM:HI}.
\vskip 9pt
\noindent{\underbar{Proof of Theorem \ref{THM:HI}.}}
\vskip 9pt
\rm\enspace
Set $y_0$ $=$ $(0,r/2)$ $\in$ $\mathbb{R}^d$.
Take the difference between \eqref{EQ:GLHF} and $\triangle w_{\lambda,l}$ $=$ $0$ on $T_l^+$,
multiplying it by $-2(y-y_0)\cdot D((v_\lambda-h_{x_0})\,-\,w_{\lambda,l})/m^2$,
integrate it on $T_l^+$ and sum up it for $l$ and $m$ to verify
\allowdisplaybreaks\begin{align}
%1
& -2\,\sum_{l=1}^{L_\lambda}\sum_{m=1}^{1/\triangle\tau_{l-1}^2}\frac 1{m^2}
\lint_{t_{l,m-1}}^{t_{l,m}}\, dt\lint_{T_l^+}
\left\langle\frac {\partial v_\lambda}{\partial t},
(y-y_0)\!\cdot\! D((v_\lambda-h_{0})\,-\,w_{\lambda,l})\right\rangle\,\, dy
\notag\\
%2
&
\,+\,(d-2)\,\sum_{l=1}^{L_\lambda}\sum_{m=1}^{1/\triangle\tau_{l-1}^2}\frac 1{m^2} 
\lint_{t_{l,m-1}}^{t_{l,m}}\, dt\lint_{T_l^+}
|D((v_\lambda-h_{x_0})\,-\,w_{\lambda,l})|^2\,\, dy\notag\\
%3
&
\,+\,\,\sum_{l=1}^{L_\lambda}\sum_{m=1}^{1/\triangle\tau_{l-1}^2}\frac 1{m^2}
\lint_{t_{l,m-1}}^{t_{l,m}}\, dt\lint_{T_l^+}
\langle a\cdot D((v_\lambda-h_{x_0})\,-\,w_{\lambda,l}),
D_d ((v_\lambda-h_{x_0})\,-\,w_{\lambda,l})\rangle\,\, dy\notag\\
%4
& \, + \, \frac {d}{2} \,
\sum_{l=1}^{L_\lambda}\sum_{m=1}^{1/\triangle\tau_{l-1}^2}\frac 1{m^2}
\lint_{t_{l,m-1}}^{t_{l,m}}\lambda^{1-\kappa}\, dt\lint_{T_l^+}
(|v_\lambda|^2\,-\,1)^2 \, dy\notag\\
%R1
& \;=\;-2\sum_{l=1}^{L_\lambda}\sum_{m=1}^{1/\triangle\tau_{l-1}^2}\frac 1{m^2}
\lint_{t_{l,m-1}}^{t_{l,m}}\,\lambda^{1-\kappa}\, dt\lint_{T_l^+}
(|v_\lambda|^2\,-\,1)\la v_\lambda, (y-y_0)\cdot\! Dw_{\lambda,l}\ra\, dy\notag\\
%R2
& \,-\, 2\sum_{l=1}^{L_\lambda}\sum_{m=1}^{1/\triangle\tau_{l-1}^2}\frac 1{m^2} 
\lint_{t_{l,m-1}}^{t_{l,m}}\,\lambda^{1-\kappa}\, dt\lint_{T_l^+} (|v_\lambda|^2\,-\,1)
\la v_\lambda,(y-y_0)\!\cdot\! Dh_{x_0}\ra\, dy\notag\\
%R3
& \,-2\,\sum_{l=1}^{L_\lambda}\sum_{m=1}^{1/\triangle\tau_{l-1}^2}
\frac 1{m^2}\lint_{t_{l,m-1}}^{t_{l,m}}\, dt\lint_{T_l^+}
\la (L\,-\,\triangle)w_{\lambda,l},(y-y_0)\!\cdot\! D((v_\lambda-h_{x_0})\,-\,w_{\lambda,l})\ra\, dy
\notag\\
%R4
& \,-\,\sum_{l=1}^{L_\lambda}\sum_{m=1}^{1/\triangle\tau_{l-1}^2}\frac 1{m^2}\,
\lint_{t_{l,m-1}}^{t_{l,m}}\, dt\lint_{\partial T_l\cap\{y_d>0\}}
(y-y_0)\cdot\nu 
|D_\nu ((v_\lambda-h_{x_0})\,-\,w_{\lambda,l})|^2 \, d\mathcal{H}_{y}^{d-1}
\notag\\
%R5
& \,+\,\sum_{l=1}^{L_\lambda}\sum_{m=1}^{1/\triangle\tau_{l-1}^2}\frac 1{m^2}\,
\lint_{t_{l,m-1}}^{t_{l,m}}\, dt\lint_{\partial T_l\cap\{y_d>0\}}
(y-y_0)\cdot\nu 
|D_\tau ((v_\lambda-h_{x_0})\,-\,w_{\lambda,l})|^2 \, d\mathcal{H}_{y}^{d-1}
\notag\\
%R6
& \,-\,2\sum_{l=1}^{L_\lambda}\sum_{m=1}^{1/\triangle\tau_{l-1}^2}\frac 1{m^2}\,
\lint_{t_{l,m-1}}^{t_{l,m}}\, dt\lint_{\partial T_l\cap\{y_d>0\}}
\langle D_\nu ((v_\lambda-h_{x_0})\,-\,w_{\lambda,l}),(y-y_0)\cdot D_\tau ((v_\lambda-h_{x_0})\,-\,w_{\lambda,l})\rangle
\, d\mathcal{H}_{y}^{d-1}\notag\\
%R7
& \,+\,2
\sum_{l=1}^{L_\lambda}\sum_{m=1}^{1/\triangle\tau_{l-1}^2}\frac 1{m^2}\,
\lint_{t_{l,m-1}}^{t_{l,m}}\, dt\lint_{\partial T_l\cap\{y_d>0\}}
\frac{y_d}{|y'|}\, \notag\\
%R8
& \qquad\qquad
\langle a\cdot D((v_\lambda-h_{x_0})\,-\,w_{\lambda,l}),(y-y_0)\cdot D((v_\lambda-h_{x_0})\,-\,w_{\lambda,l})\rangle
\, d\mathcal{H}_{y}^{d-1}\notag\\
%R9
& \,-\,
\sum_{l=1}^{L_\lambda}\sum_{m=1}^{1/\triangle\tau_{l-1}^2}\frac 1{m^2}\,
\lint_{t_{l,m-1}}^{t_{l,m}}\, dt\lint_{\partial T_l\cap\{y_d>0\}}
\frac{y_d}{|y'|}\,a\cdot y |D((v_\lambda-h_{x_0})\,-\,w_{\lambda,l})|^2 \, d\mathcal{H}_{y}^{d-1}\notag\\
%R10
& \,-\,\sum_{l=1}^{L_\lambda}\sum_{m=1}^{1/\triangle\tau_{l-1}^2}\frac 1{m^2}\,
\left.
\lint_{t_{l,m-1}}^{t_{l,m}}\, dt
\lint_{B_\rho\cap\{y_d=0\}}(1+|a|^2)\bigl(\frac\rho 2+a\cdot y\bigr)
|D_d ((v_\lambda-h_{x_0})\,-\,w_{\lambda,l})|^2\, d\mathcal{H}_{y}^{d-1}
\right|_{\rho=r_{l-1}}^{r_l}\notag\\
%R11
&
\,+\,\sum_{l=1}^{L_\lambda}\sum_{m=1}^{1/\triangle\tau_{l-1}^2}\frac {1}{2{m^2}}
\lint_{t_{l,m-1}}^{t_{l,m}}\,{\lambda^{1-\kappa}}\, dt\lint_{\partial T_l^+}
\bigl(\nu\cdot(y-y_0)\,-\,\frac{a\!\cdot\! y\, y_d}{|y'|}\bigr) 
(|v_\lambda|^2\,-\,1)^2\, d\mathcal{H}_y^{d-1}\notag\\
%R11
&
\;=\;\;(\mathrm{\bigroman{1}})\,+\,(\mathrm{\bigroman{2}})\,+\, \,(\mathrm{\bigroman{3}})
\,+\,\cdots\,+\,(\mathrm{\bigroman{10}})\label{EQ:RP-0}\end{align}
with \begin{math}y\,=\,(y',y_d)\,\in\,\mathbb{R}^d,\end{math}
\begin{math}
a\!\cdot\!y \; = \; \sum_{i=1}^{d-1} a_i y_i, \; a\!\cdot\!D \; = \; \sum_{i=1}^{d-1} a_i D_i
\end{math} and
\begin{math}a\!\cdot\! D_\tau\;=\;\sum_{i=1}^{d-1} a_i D_{\tau_i}.\end{math}
Here we used 
\begin{math}
\nabla_i a_j\,-\,\nabla_j a_i\;=\;(\nabla_i\nabla_j\,-\,\nabla_j\nabla_i) \phi\;=\;0.
\end{math}
\par
From now on, we shall estimate each term of the right-hand side in \eqref{EQ:RP-0}.
Primarily we estimate ({\bigroman{1}}): By recalling Remark \ref{REM:DensityOfGLQE}, we have
\begin{align}
%1
(\mathrm{\bigroman{1}}) &\;\le\;\epsilon_0^2\lint_{Q_r} \lambda^{1-\kappa} 
( 1 \, - \, |u_\lambda|^2 )\, dz\notag\\
%2
& \,+\, C(\epsilon_0)\lint_0^{r^2} \lambda^{1-\kappa}\, dt
\lint_{B_{(1-\epsilon_0^4)r}} ( 1 \, - \, |v_\lambda|^2 )
\fint\limits_{\partial B_{r}} |v_\lambda-h_{x_0}|^2 \, d\mathcal{H}_{y}^{d-1}
\notag\\
%3
&
\,+\,\frac 1{\epsilon_0^2}\sum_{l=2}^{L_\lambda-2}\sum_{m=1}^{1/\triangle\tau_{l-1}^2} \frac 1{m^2}
\lint_{t_{l,m-1}}^{t_{l,m}} \lambda^{1-\kappa} \, dt\, \lint_{T_l^+}
(1\,-\,|u_\lambda|^2) \, |(y-y_0)\!\cdot\! Dw_{\lambda,l}|^2 \, dy\notag\\
%4
& \,+\, o(1)\qquad (\lambda\nearrow\infty).\label{INEQ:R-1pre}\end{align}
\par We estimate the third term above called (\romannumeral 3):
Use a polar coordinate \begin{math}x\;=\;(|x|,x/|x|)\;=\;(\rho,\omega_{d-1})\end{math} and
set $S_\tau^{d-1}(\omega_{d-1})$ being a geodesic ball on $\mathbb{S}^{d-1}$
with the radius $\tau$ and the center $\omega_{d-1}$.
For any $l$ $(l=2,3,\ldots,L_\lambda-1)$, 
we choose a collection of balls $\{S_{\triangle{\tau_l}} (\omega_{d-1}^j)\}$ $(j\,\in\, J_l)$ 
for some positive integer $J_l$ with
\allowdisplaybreaks\begin{align}
%1
&\mathbb{S}^{d-1}\,\subset\,\bigcup_{j\in J_l}S_{\triangle{\tau_l}} (\omega_{d-1}^j),\notag\\
%2
& \max_{j\,\in\, J_l}
\mathrm{Card}\{j\,\in\, J_l\,;\,
S_{\triangle{\tau_l}}(\omega_{d-1}^j)\cap S_{\triangle{\tau_l}}(\omega_{d-1}^k)\;\ne\;\emptyset\}\;<\; C \;<\;\infty
\notag\end{align}
for each $k$ $=$ $1,2,\ldots,J_l$,
where $C$ is a positive constant independent of $k$ and $l$.
Since 
\allowdisplaybreaks\begin{align}
%1
\sup_{[r_{l-1},r_l)\times S_{\triangle{r_l}}(\omega_{d-1}^j)}&
|\nabla_\nu w_{\lambda,l}|^2\;\le\;
C \triangle r_l^2\sup_{[\widetilde{r}_{l-1},\widetilde{r}_{l})\times S_{\triangle\tau_l}(\omega_{d-1}^j)}
|\nabla_\nu f_{\lambda}|^2\notag\\
%2
& \;\le\;\frac C{\triangle{r_l}^{d-1}}
\fint\limits_{t_{l,m-1}}^{t_{l,m}} \, dt 
\lint_{[\widetilde{r}_{l-2},\widetilde{r}_{l+1})\times S_{2\triangle\tau_l}(\omega_{d-1}^j)}
|\nabla_\nu f_\lambda|^2 \, dy \label{INEQ:Key}\\
%3
& \,+\,\frac {C \triangle{r_l}^{3}}{\triangle{r_l}^{d}}
\lint\limits_{t_{l,m-1}}^{t_{l,m}} \, dt
\lint_{[\widetilde{r}_{l-2},\widetilde{r}_{l+1})\times S_{2\triangle\tau_l}(\omega_{d-1}^j)}
|\nabla_\nu \nabla_t f_\lambda|^2 \, dy
\notag\end{align}
holds, employing Remark \ref{REM:DensityOfGLQE}, 
the third term $(\mathrm{\romannumeral 3})$ is estimated as follows:
\allowdisplaybreaks\begin{align}
%1
& \epsilon_0^2 (\mathrm{\romannumeral 3})\;\le\; r^2
\sum_{l=2}^{L_\lambda-2}\sum_{m=1}^{1/\triangle\tau_{l-1}^2} \frac 1{m^2}\sum_{j\in J_l} 
\lint_{t_{l,m-1}}^{t_{l,m}}\, {\lambda^{1-\kappa}}\, dt \kern-15pt
\lint_{[{r}_{l-1},{r}_{l})\times S_{\triangle\tau_l}(\omega_{d-1}^j)}\kern-15pt
(1\,-\,|v_\lambda |^2 )\, dy \,
\sup_{[\widetilde{r}_{l-1},\widetilde{r}_{l})\times S_{\triangle\tau_l}(\omega_{d-1}^j)}|\nabla w_{\lambda,l}|^2
\notag\\
%2
& \;\le\; Cr^2
\sum_{l=2}^{L_\lambda-2}{\triangle r_l}\sum_{m=1}^{1/\triangle\tau_{l-1}^2}\frac 1{m^2}
\sum_{j\in J_l}\frac 1{\triangle r_l^d} 
\lint_{{t}_{l,m-1}}^{{t}_{l,m}}\lambda^{1-\kappa}\, dt \, 
\lint_{[{r}_{l-1},{r}_{l})\times S_{\triangle\tau_l}(\omega_{d-1}^j)} (1\,-\, |v_\lambda |^2 )\, dy
\notag\\
%3
&\qquad\qquad\fint_{{t}_{l,m-1}}^{{t}_{l,m}} \, dt
\lint_{[\widetilde{r}_{l-2},\widetilde{r}_{l+1})\times\mathbb{S}_{2\triangle\tau_l}^{d-1}(\omega_{d-1}^j)}
|\nabla_\nu f_\lambda |^2 \, dy\notag\\
%4
& \,+\, Cr^2
\sum_{l=2}^{L_\lambda-2}{\triangle r_l}\sum_{m=1}^{1/\triangle\tau_{l-1}^2}\frac 1{m^2}
\sum_{j\in J_l}\frac 1{\triangle r_l^d}
\lint_{{t}_{l,m-1}}^{{t}_{l,m}}\lambda^{1-\kappa}\, dt \,
\lint_{[{r}_{l-1},{r}_{l})\times S_{\triangle\tau_l}(\omega_{d-1}^j)} (1\,-\, |v_\lambda |^2 )\, dy
\notag\\
%5
& \qquad\qquad\fint_{{t}_{l,m-1}}^{{t}_{l,m}} \, dt
\lint_{\{\widetilde{r}_{l+1}\}\times\mathbb{S}_{2\triangle\tau_l}^{d-1}(\omega_{d-1}^j)}
|\nabla_\tau f_\lambda|^2\, d\mathcal{H}_y^{d-1}\notag\\
%6
& \,+\, Cr^2
\sum_{l=2}^{L_\lambda-2}\triangle r_l^3\sum_{m=1}^{1/\triangle\tau_{l-1}^2}\frac 1{m^2}
\sum_{j\in J_l}\frac 1{\triangle r_l^d}
\lint_{{t}_{l,m-1}}^{{t}_{l,m}}\lambda^{1-\kappa}\, dt \,
\lint_{[{r}_{l-1},{r}_{l})\times S_{\triangle\tau_l}(\omega_{d-1}^j)} (1\,-\, |v_\lambda |^2 )\, dy\notag\\
%7
& \qquad\qquad \lint_{{t}_{l,m-1}}^{{t}_{l,m}} \, dt
\lint_{[\widetilde{r}_{l-2},\widetilde{r}_{l+1})\times\mathbb{S}_{2\triangle\tau_l}^{d-1}(\omega_{d-1}^j)}
|\nabla_\nu \nabla_t f_\lambda |^2 \, dy\notag\\
%8
& \,+\, Cr^2
\sum_{l=2}^{L_\lambda-2}\triangle r_l^3\sum_{m=1}^{\max (1,1/\triangle\tau_{l-1}^3)}\frac 1{m^2}
\sum_{j\in J_l}\frac 1{\triangle r_l^d}
\lint_{{t}_{l,m-1}}^{{t}_{l,m}}\lambda^{1-\kappa}\, dt \,
\lint_{\widetilde{T}_l} (1\,-\, |v_\lambda |^2 )\, dy\notag\\
%9
&\qquad\qquad \lint_{{t}_{l,m-1}}^{{t}_{l,m}}\, dt
\lint_{\{\widetilde{r}_{l+1}\}\times\mathbb{S}_{2\triangle\tau_l}^{d-1}(\omega_{d-1}^j)}
|\nabla_\tau \nabla_t f_\lambda|^2 \, d\mathcal{H}_y^{d-1}\notag\\
%10
&\;\le\;Cr^2\sum_{l=1}^{L_\lambda-1}{\triangle r_l}\sum_{m=1}^{1/\triangle\tau_{l-1}^2}\frac 1{m^2}
\fint_{{t}_{l,m-1}}^{{t}_{l,m}} \, dt
\lint_{\widetilde{T}_{l}} |\nabla_\nu f_\lambda |^2 \, dy\notag\\
%11
&\,+\,Cr^2\sum_{l=3}^{L_\lambda-1}\triangle r_l\sum_{m=1}^{1/\triangle\tau_{l-1}^2}\frac 1{m^2}
\fint_{{t}_{l,m-1}}^{{t}_{l,m}}\, dt
\lint_{\partial\widetilde{B}_{l}}|\nabla_\tau f_\lambda|^2\, dy\notag\\
%12
&\,+\,Cr^2
\sum_{l=1}^{L_\lambda-1}\triangle r_l^3\sum_{m=1}^{1/\triangle\tau_{l-1}^2}\frac 1{m^2}
\lint_{{t}_{l,m-1}}^{{t}_{l,m}}\, dt \,
\lint_{\widetilde{T}_{l}}|\nabla_\nu\nabla_t f_\lambda |^2\, dy\notag\\
%13
&\,+\,Cr^2
\sum_{l=3}^{L_\lambda-1}\triangle r_l\sum_{m=1}^{1/\triangle\tau_{l-1}^2}\frac 1{m^2}
\lint_{{t}_{l,m-1}}^{{t}_{l,m}}\, dt \,
\lint_{\widetilde{T}_{l}}|\nabla_\tau\nabla_t f_\lambda|^2\, dy\notag\\
%14
&\;\le\;C\epsilon_0^4 r^3
\sup_{t\in [0,r^2)}\lint_{\partial B_{r}}|\nabla_\tau (v_\lambda\,-\,h_{x_0})|^2\, d\mathcal{H}_y^{d-1}
\,+\,C\epsilon_0^4 r^3
\lint_{0}^{r^2}\, dt\, \lint_{\partial B_{r}} 
\left|\frac{\partial v_\lambda}{\partial t}\right|^2 \, d\mathcal{H}_y^{d-1}.
\notag\end{align}
Thus by applying Lemma \ref{LEM:GLQE} to the second term, (\bigroman{1}) is estimated as follows:
\allowdisplaybreaks\begin{align}
%1
& (\mathrm{\bigroman{1}})\;\le\; \epsilon_0^2\lint_{Q_r}\lambda^{1-\kappa}(1\,-\,|v_\lambda |^2)\, dz
\notag\\
%2
&\,+\,C\epsilon_0^2 r^3\sup_{t\in[0,r^2)}\lint_{\partial B_r}|\nabla (v_\lambda\,-\,h_{0})|^2\, d\mathcal{H}_y^{d-1}
\,+\,C\epsilon_0^2 r^3\lint_0^{r^2}\, dt\lint_{\partial B_r} 
\left|\frac{\partial v_\lambda}{\partial t}\right|^2 \, dy \notag\\
%3
& \,+\,C(\epsilon_0) r\sup_{t\in[0,r^2)}\lint_{\partial B_r} |v_\lambda\,-\, h_{x_0}|^2\, d\mathcal{H}_y^{d-1}.
\label{INEQ:R-1}\end{align}
\par
Next we asses (\bigroman{2}): 
By virtue of D.Gilbarg and N.S.Trudinger~\cite[Theorem 14.2]{gilbarg-trudinger}, 
we infer
\allowdisplaybreaks\begin{align}
%1
(\mathrm{\bigroman{2}}) & \;\le\; Cr \lint_{Q_r} e_\lambda \,dz
\,+\, o(1)\qquad (\lambda\nearrow\infty). \label{INEQ:R-2}\end{align}
Next we estimate the third term:
\allowdisplaybreaks\begin{align}
%1
(\mathrm{\bigroman{3}}) & \; =\; 2\lint_0^{r^2}\, dt\lint_{T_1^+}
\langle (L-\triangle) w_{\lambda,1,1},(y-y_0)\cdot D ((v_\lambda-h_{x_0})\,-\,w_{\lambda,1,1})\rangle\, dy
\notag\\
%2
&\,+\,2\sum_{l=2}^{L_\lambda-1}\sum_{m=1}^{1/\triangle\tau_{l-1}^2}\frac 1{m^2}
\lint_{t_{l,m-1}}^{t_{l,m}}\, dt \lint_{T_l^+} 
\langle	(L-\triangle) w_{\lambda,l}, (y-y_0)\cdot D ((v_\lambda-h_{x_0})\,-\,w_{\lambda,l})\rangle\, dy
\notag\\
%3
&\;\le\;C\epsilon_0^2
\lint_0^{r^2}\, dt\lint_{B_1^+}|\nabla((v_\lambda-h_{x_0})\,-\,w_{\lambda,l})|^2\, dy
\,+\,\frac{r^2}{2\epsilon_0^2}\lint_0^{r^2}\, dt\lint_{T_1^+} |(L-\triangle)w_{\lambda,1,1}|^2\, dy
\notag\\
%4
& \,+\,\frac{r^2}{\epsilon_0^2} \sum_{l=2}^{L_\lambda-1}\sum_{m=1}^{1/\triangle\tau_{l-1}^2}\frac 1{m^2}
\lint_{t_{l,m-1}}^{t_{l,m}}\, dt \lint_{T_l^+}
|(L-\triangle)w_{\lambda,l}|^2 \, dy. \label{INEQ:R-3-1}\end{align}
We estimate the final term of \eqref{INEQ:R-3-1}: 
A direct computation from the expression \eqref{EQ:Formula} reads
\allowdisplaybreaks\begin{align}
%1
& r^2 \sum_{l=2}^{L_\lambda-2}\sum_{m=1}^{1/\triangle\tau_{l-1}^2}\frac 1{m^2}
\lint_{t_{l,m-1}}^{t_{l,m}}\, dt \lint_{T_l^+}|(L-\triangle)w_{\lambda,l}|^2 \, dy\notag\\
%2
& \;\le\; C\sum_{l=2}^{L_\lambda-2}\triangle{r}_l \sum_{m=1}^{1/\triangle\tau_{l-1}^2}\frac 1{m^2}
\lint_{t_{l,m-1}}^{t_{l,m}}\, dt \lint_{\partial\widetilde{B}_l^+}
(|\nabla_\tau f_\lambda|^2\,+\,\triangle{r}_l^2 |\nabla_\nu\nabla_\tau f_\lambda|^2)
\, d\mathcal{H}_y^{d-1}\notag\\
%3
&\;\le\;C\sum_{l=2}^{L_\lambda-2}\triangle{r}_l^3
\sup_{[0,r^2)} \lint_{{T}_{l+1}^+}|\nabla_\tau f_\lambda|^2
\,+\,C\sum_{l=2}^{L_\lambda-2}\triangle{r}_l^3
\sup_{[0,r^2)}\lint_{\widetilde{T}_{l+1}^+}|\nabla_\nu\nabla_\tau f_\lambda|^2 \, dy\notag\\
%4
&\;\le\;C\epsilon_0^4 r^3\sup_{[0,r^2)}\lint_{\partial B_{r}^+}
|\nabla_\tau (v_\lambda-h_{x_0})|^2\, d\mathcal{H}_y^{d-1}.\end{align}
Hence we conclude
\allowdisplaybreaks\begin{align}
%1
(\mathrm{\bigroman{3}}) & \;\le\; C\epsilon_0^2\lint_{Q_r^+} |\nabla (v_\lambda-h_{x_0})|^2 \, dz
\,+\,\frac{C(\epsilon_0)}r \lint_0^{r^2}\, dt \lint_{\partial B_r^+} |v_\lambda\,-\,h_{x_0}|^2 \, d\mathcal{H}_y^{d-1}
\notag\\
%2
&\,+\,C\epsilon_0^2 r^3 \sup_{[0,r^2)}\lint_{\partial B_r^+}
|\nabla (v_\lambda\,-\,h_{x_0})|^2 \, d\mathcal{H}_{y}^{d-1}.
\label{INEQ:R-3}\end{align}
\par 
Successively we estimate (\bigroman{4}): 
\allowdisplaybreaks\begin{align}
%1
& (\mathrm{\bigroman{4}}) \;=\;-
\sum_{l=1}^{L_\lambda}\sum_{m=1}^{1/\triangle\tau_{l-1}^2}\frac 1{m^2}
\lint_{t_{l,m-1}}^{t_{l,m}}\, dt\lint_{\partial B_{l}^+} (y-y_0)\cdot\nu
|D_{\nu}((v_\lambda\,-\,h_{x_0})\,-\,w_{\lambda,l})|^2\,d\mathcal{H}_y^{d-1}\notag\\
%2
&\,+\,\sum_{l=1}^{L_\lambda-1}\sum_{m=1}^{1/\triangle\tau_{l}^2}\frac 1{m^2}
\lint_{t_{l+1,m-1}}^{t_{l+1,m}}\, dt\lint_{\partial B_{l}^+} (y-y_0)\cdot\nu
|D_{\nu}((v_\lambda\,-\,h_{x_0})\,-\,w_{\lambda,1+1})|^2\,d\mathcal{H}_y^{d-1}\notag\\
%3
& \;\le\; -\frac r2\sum_{m=1}^{1/\triangle\tau_{L_\lambda-1}^2}\frac 1{m^2}
\lint_{t_{L_\lambda,m-1}}^{t_{L_\lambda,m}}\, dt\lint_{\partial B_{L_\lambda}^+}
|D_{\nu}((v_\lambda\,-\,h_{x_0})\,-\,w_{\lambda,L_\lambda})|^2\,d\mathcal{H}_y^{d-1}\notag\\
%4
&\,+\, 2\lint_{0}^{r^2}\, dt\lint_{\partial B_{1}^+} (y-y_0)\cdot\nu
\langle D_\nu(w_{\lambda,1}\,-\,w_{\lambda,2}),
D_\nu( (v_\lambda-h_{x_0})\,-\, (w_{\lambda,1}\,+\,w_{\lambda,2})/2)\rangle\,d\mathcal{H}_y^{d-1}\notag\\
%5
&\,+\,Cr\sum_{l=2}^{L_\lambda-1}\sum_{m=1}^{1/\triangle\tau_{l-1}^2}\frac 1{m^2}
\lint_{t_{l,m-1}}^{t_{l,m}}\, dt\lint_{\partial B_{l-1}^+} 
|D_{\nu}((v_\lambda\,-\,h_{x_0})\,-\,w_{\lambda,l})|^2\,d\mathcal{H}_y^{d-1}\notag\\
\notag\\
%6
&\;\le\;-\frac r2 
\sum_{m=1}^{1/\triangle\tau_{L_\lambda-1}^2}\frac 1{m^2}
\lint_{t_{L_\lambda,m-1}}^{t_{L_\lambda,m}}\, dt\lint_{\partial B_{L_\lambda}^+}
|D_{\nu}((v_\lambda\,-\,h_{x_0})\,-\,w_{\lambda,L_\lambda})|^2\,d\mathcal{H}_y^{d-1}\notag\\
%7
&\,+\,\epsilon_0^2\lint_0^{r^2}\, dt\lint_{\partial B_1^+}|\nabla (v_\lambda-h_{x_0})|^2 \, dx
\,+\,\frac{C(\epsilon_0)}r \lint_0^{r^2}\, dt\lint_{\partial B_r^+}
|v_\lambda-h_{x_0}|^2\, d\mathcal{H}_y^{d-1}\notag\\
%8
&\,+\,Cr\sum_{l=1}^{L_\lambda-1}\triangle r_{l}
\sum_{m=1}^{1/\triangle\tau_{l-1}^2}\frac 1{m^2}
\sup_{t\in [0,r^2)}\lint_{T_l^+} |\nabla (v_\lambda\,-\,h_{x_0})|^2\,dy \notag\\
%9
& \,+\,Cr\sum_{l=2}^{L_\lambda-1}\triangle r_{l}^2
\sum_{m=1}^{1/\triangle\tau_{l-1}^2}\frac 1{m^2}
\sup_{t\in [0,r^2)}\lint_{\partial B_{r_{L_\lambda}}^+}
|\nabla_\tau (v_\lambda\,-\,h_{x_0})|^2\,d\mathcal{H}_y^{d-1} \notag\\
%10
& \,+\,Cr\sum_{l=2}^{L_\lambda-1}\triangle r_{l}^2
\sum_{m=1}^{1/\triangle\tau_{l-1}^2}\frac 1{m^2}
\sup_{t\in [0,r^2)}\lint_{\widetilde{T}_l^+} |\nabla_\nu f_\lambda|^2\,dy\notag\\
%11
&\;\le\;-\frac r2
\sum_{m=1}^{1/\triangle\tau_{L_\lambda-1}^2}\frac 1{m^2}
\lint_{t_{L_\lambda,m-1}}^{t_{L_\lambda,m}}\, dt\lint_{\partial B_{L_\lambda}^+}
|D_{\nu}((v_\lambda\,-\,h_{x_0})\,-\,w_{\lambda,L_\lambda})|^2\,d\mathcal{H}_y^{d-1}\notag\\
%12
&\,+\,\epsilon_0^2\lint_{Q_r^+}|\nabla (v_\lambda\,-\,h_{x_0})|^2\,dz\notag\\
%13
&\,+\,C\epsilon_0^4 r^2\sup_{t\in [0,r^2)}\lint_{B_{r}^+}|\nabla (v_\lambda\,-\,h_{x_0})|^2 \, dy
\,+\,C\epsilon_0^4 r^3 \sup_{t\in [0,r^2)}\lint_{\partial B_{r}^+}
|\nabla_{\tau} (v_\lambda\,-\,h_{x_0})|^2 \, d\mathcal{H}_y^{d-1}\notag\\
%14
& \,+\, \frac {C(\epsilon_0)}r \lint_{0}^{r^2}\, dt\lint_{\partial B_{r}^+}
|v_\lambda\,-\, h_{x_0}|^2 \, d\mathcal{H}_y^{d-1}. \label{INEQ:R-4}\end{align}
\par
We have a similar estimate for $(\mathrm{\bigroman{5}})$, $\ldots$, $(\mathrm{\bigroman{8}})$.
So $(\mathrm{\bigroman{4}})$ $+\cdots+$ $(\mathrm{\bigroman{8}})$ are majorized above as follows:
On account of a smallness of $R$ and a continuity of $a$ $=$ $\nabla\phi$,
there exists a positive number $\theta_0$ less than $1$ such that
\allowdisplaybreaks\begin{align}
%1
& (\mathrm{\bigroman{4}})\,+\,(\mathrm{\bigroman{5}})\,+\,\cdots\,+\,(\mathrm{\bigroman{8}})\notag\\
%2
&\;\le\;-\theta_0\sum_{m=1}^{1/\triangle\tau_{L_\lambda-1}^2}\frac 1{m^2}
\lint_{t_{L_\lambda,m-1}}^{t_{L_\lambda,m}}\, dt\lint_{\partial B_{l}^+}
|D_{\nu}((v_\lambda\,-\,h_{x_0})\,-\,w_{\lambda,L_\lambda})|^2\,d\mathcal{H}_y^{d-1}\notag\\
%3
&\,+\,C\epsilon_0^4 r^2\sup_{t\in [0,r^2)}\lint_{B_{r}^+}|\nabla (v_\lambda\,-\,h_{x_0})|^2 \, dy
\,+\,
C\epsilon_0^4 r^3 \sup_{t\in [0,r^2)}\lint_{\partial B_{r}^+}
|\nabla_{\tau} (v_\lambda\,-\,h_{x_0})|^2 \, d\mathcal{H}_y^{d-1}\notag\\
%4
& \,+\, \frac {C(\epsilon_0)}r \lint_{0}^{r^2}\, dt\lint_{\partial B_{r}^+}
|v_\lambda\,-\, h_{x_0}|^2 \, d\mathcal{H}_y^{d-1}. \label{INEQ:R-4+5+6+7+8+9}\end{align}
Finally we know that $(\mathrm{\bigroman{10}})$ becomes
\allowdisplaybreaks\begin{align}
%1
& (\mathrm{\bigroman{10}})\;=\;\lint_0^{r^2}\frac{r \lambda^{1-\kappa}}2 \, dt
\lint_{\partial B_r^+} \bigl(1\,-\,\frac{a\cdot\nu y_d}{|y'|}\bigr)
(|u_\lambda|^2\,-\,1)^2\, d\mathcal{H}_y^{d-1}.
\label{INEQ:R-10}\end{align}
\par
On the other hand, we also have the following estimate for the left-hand side in \eqref{EQ:RP-0},
which is called $(L)$:
\allowdisplaybreaks\begin{align}
%1
\mathrm{(L)}& \,\ge\;\frac{d-2}4\lint_{B_{(1-\epsilon_0^4)r}^+} \mathbf{e}_\lambda \, dy
\,-\,C\,r^2\lint_{B_r^+}\left|\frac {\partial v_\lambda}{\partial t}\right|^2 \,  \, dy
\label{INEQ:L}\\
%2
&
\,-\,\frac{C (\epsilon_0)}{r}\lint_{\partial B_r^+}|v_\lambda-h_{x_0}|^2\, d\mathcal{H}_y^{d-1}
\,-\,C\lint_{B_r^+}\vert\nabla h_{x_0}\vert^2\, dy.\notag\end{align}
\par
A substitution of \eqref{INEQ:R-1}, \eqref{INEQ:R-2}, \eqref{INEQ:R-3},\eqref{INEQ:R-4+5+6+7+8+9} 
and \eqref{INEQ:R-10} for \eqref{EQ:RP-0} verifies 
\allowdisplaybreaks\begin{align}
%1
\lint_{Q_{r/2}^+} & \mathbf{e}_\lambda\, dz
\;\le\; C\epsilon_0^2\lint_{P_{r}^+}\mathbf{e}_\lambda \, dz
\,+\,\frac {C r^2}{\epsilon_0^2}\lint_{P_{r}^+} \biggl|\frac{\partial v_\lambda}{\partial t}\biggr|^2\, dz
\,+\,{C r^3}\lint_{0}^{r^2}\, dt\lint_{\partial B_r^+}
\left\vert\frac{\partial v_\lambda}{\partial t}\right\vert^2 \, d\mathcal{H}_y^{d-1}
\notag\\
%2
&\,+\,{C \epsilon_0^2 r^3}\sup_{[0,r^2)}\lint_{\partial B_r^+}
\vert\nabla(v_\lambda\,-\,h_{x_0})\vert^2\, d\mathcal{H}_y^{d-1}
\,+\,\frac {C(\epsilon_0)}r 
\lint_0^{r^2}\, dt\lint\limits_{\partial B_{r}^+}
\vert v_\lambda\,-\,h_{x_0}\vert^2\, d\mathcal{H}_y^{d-1}\notag\\  
%3
&
\,+\, Cr\lint_{-r^2}^{r^2}{\lambda^{1-\kappa}}\, dt  
\,\lint_{\partial B_{r}^+} ( |v_\lambda|^2 \, - \, 1 )^2 \, d\mathcal{H}_y^{d-1}
\,+\,C\lint_{P_{r}^+}\vert\nabla h_{x_0}\vert^2\, dz\,+\, o(1)
\quad (\lambda\,\nearrow\,\infty), \label{INEQ:RP-2}\end{align} 
where we used a smallness of $R$ to estimate (\bigroman{2}).
\par
A way of choosing $r$ deduces
\allowdisplaybreaks\begin{align}
%1
\lint_{Q_{R}^+} & \mathbf{e}_\lambda \, dz
\;\le\; C\epsilon_0^2\lint_{P_{2R}^+} \mathbf{e}_\lambda \, dz
\,+\,\frac {C R^2}{\epsilon_0^2}\lint_{P_{2R}^+}\,\Bigl|\frac{\partial v_\lambda}{\partial t}\Bigr|^2\, dz
\label{INEQ:RP-4}\\
%2
&\,+\,\frac {C(\epsilon_0)}{R^2}\lint_{P_{2R}^+}|v_\lambda-h_{x_0}|^2\, dz
\,+\,C\lint_{P_{2R}^+}\vert\nabla h_{x_0}\vert^2\, dz\,+\, o(1)
\quad (\lambda\,\nearrow\,\infty). \notag\end{align}
\par The rest of our proof is as same as the one by K.Horihata~\cite{horihata}.

%% file: bhf_bhf.tex
%
%#! platex bhf
%
%%%%%%%%%%%%%%%%%%%%%%%%%%%%%%%%%%%%%%%%%%%%%%%%%%%%%%%%%%%%%%%%%%%%%%%%%%%
\setcounter{chapternumber}{2}\setcounter{equation}{0}
\renewcommand{\theequation}%
{\thechapternumber.\arabic{equation}}
\setcounter{chapternumber}{3}\setcounter{equation}{0}
\renewcommand{\theequation}%
{\thechapternumber.\arabic{equation}}
\section{\enspace WHHF}
\subsection{\enspace Partial Regularity}
This chapter studies a partial boundary regularity on WHHF.
We can regard all theorems, lemmas and corollary as the boundary version to
chapter 3 in K.Horihata\cite{horihata}.
We commence with reviewing the following existence theorem:
\begin{Thm}{\rm{(Existence).}}\label{THM:Existence-HF}
The GLHF converges to a WHHF in $L^2 (Q(T))$ 
as $\lambda \nearrow \infty$ \rm{(}\it{modulo a sub-sequence of $\lambda$}\rm{)}.
\end{Thm}
\begin{Def}{}\label{REM:Measure}
Let $\{u_{\lambda (\nu )}\}$ $( \nu \, = \, 1,2,\ldots )$
be the sequence selected above and set $\mathbf{e}_{\lambda (\nu)}$ as
the Ginzburg-Landau energy density
$|\nabla u_{\lambda (\nu)}|^2/2$ $+$
$\lambda (\nu)^{1-\kappa }$ $(|u_{\lambda (\nu)}|^2\,-\,1)^2/4$.
We then denotes $\overline{\mathcal{M}}$ by
\begin{equation*}
\overline{\mathcal{M}} (P_R (z_0)\cap Q(T)) \,= \, 
\limsup_{\lambda(\nu)\nearrow\infty}\frac 1{R^d}\lint_{P_R(z_0)\cap Q(T)}
\mathbf{e}_{\lambda (\nu)}\,dz\end{equation*}
for an arbitrary parabolic cylinder $P_R (z_0)$.
\end{Def}
\begin{Lem}{\rm{(Measured Hybrid Inequality).}}\label{LEM:Measure-Hybrid}
Set \begin{math}a\;=\;\fint\limits_{B_{2R}\cap\Omega} u_0\,dx\end{math} and
assume that a sequence of GLHF $\{v_{\lambda (\nu)}\}$ 
$(\nu\, =\, 1,2,\ldots )$, converges weakly and weakly-$*$ 
in $H^{1,2} (0,T\, ;\, L^2(\Omega\,;\,\mathbb{R}^{D+1}))$ and 
$L^\infty (0,T;H^{1,2}(\Omega\,;\,\mathbb{R}^{D+1}))$
to a WHHF $u$ $\,\in\,$ $V(Q(T)\,;\,\mathbb{S}^D)$ as $\lambda (\nu)\nearrow\infty$ respectively.
Then take the pass to the limit $\lambda (\nu)\nearrow\infty$ in 
Theorem \ref{THM:HI} to infer the following{\rm{:}}
For any positive number $\epsilon_0$,
there exists a positive constant $C (\epsilon_0)$
satisfying $C (\epsilon_0)$ $\nearrow \infty$ as $\epsilon_0\searrow 0$
and depending only on $d$ and $\epsilon_0$
such that the inequality
\allowdisplaybreaks\begin{align}
%1
\overline{\mathcal{M}} & \bigl(P_R (z_0)\cap Q(T)\bigr) 
\;\le\;\epsilon_0 \overline{\mathcal{M}}\bigl(P_{4R}(z_0)\cap Q(T)\bigr)\notag\\
%2
&\,+\,C(\epsilon_0)\fint\limits_{P_{4R}(z_0)\cap Q(T)}|u\,-\,h_{x_0}|^2\,dz
\,+\, C\fint\limits_{P_{4R}\cap Q(T)} |\nabla u_0|^2 \, dz
\label{INEQ:Measure-Hybrid}\end{align}
holds for any parabolic cylinder $P_{4R} (z_0)$.\end{Lem}
Likewise L.Simon~\cite[Lemma 2, p.31]{simon95},
we can assert the reverse Poincar\'e inequality.
\begin{Cor}{\rm{(Reverse Poincar\'e Inequality).}}\label{COR:RPI}
The inequality \eqref{INEQ:Measure-Hybrid} implies
\allowdisplaybreaks\begin{align}
%1  
\overline{\mathcal{M}} (P_{R}(z_0)\cap Q(T)) &
\;\le\;C\fint\limits_{P_{2R} (z_0)\cap Q(T)} |u\,-\,h_0|^2\,dz
\notag\\
%2
& \,+\,C\fint\limits_{P_{2R}\cap{Q}(T)} |\nabla u_0|^2\, dz
\label{INEQ:RPI}\end{align}
whenever $P_{2R} (z_0)$ is an arbitrary parabolic cylinder.
\end{Cor}
By combining Corollary \ref{COR:RPI} with Sobolev imbedding theorem and Poincar\'{e} inequality,
following the proof of Theorem 2.1 in M.Giaquinta and M.Struwe~\cite{giaquinta-struwe},
we can describe the following lemma.
\begin{Lem}{}\label{LEM:RH}
Assume \begin{math}u_0\in H^{1,2q_0}(\Omega)\end{math} 
for a positive number $q_0$ greater than $1$. 
Then the WHHF $u$ belongs to $L^{2q_0}$ $((0,T)\,;\,H^{1,2q_0}(\Omega))$
with
\vskip 15pt\begin{equation}\begin{split}
%1
\biggl(\fint_{P_{R}(z_0)\cap Q(T)} |\nabla u|^{2q_0}\, dz\biggr)^{1/2q_0}
&\;\le\; C\biggl(\fint_{P_{2R}(z_0)\cap Q(T)}|\nabla u|^2\, dz\biggr)^{1/2}
\notag\\
%2
&\quad \,+\,C\biggl(\fint_{P_{2R}(z_0)\cap Q(T)}|\nabla u_0|^{2q_0}\, dz\biggr)^{1/2q_0} 
\end{split}\label{INEQ:RH}\end{equation}
for any parabolic cylinder $P_{2R}$
where $C$ is a constant depending only on $d$, $D$ and $u_0$.
\end{Lem}
If we apply \eqref{INEQ:Mon} in Theorem \ref{THM:Mon} to
the proof on a partial regularity result by Y.Chen and F.H.Lin~\cite[Theorem 2.1]{chen-lin},
then we can claim
\begin{Thm}{}\label{THM:PR}
For any positive number $\epsilon$, set
\begin{align}
%1
& \mathbf{sing}(\epsilon)\,=\,\bigcap_{R > 0}
\{z_0\,\in\,\overline{Q(T)}
\,;\,\overline{\mathcal{M}}\!\bigl(P_R (z_0)\cap Q(T)\bigr)\,\ge\,\epsilon \,\},
\label{EQ:Sing}\\
&
\mathbf{reg}(\epsilon)\;=\;\overline{Q(T)}\setminus\mathbf{sing}(\epsilon).
\label{EQ:Reg}\end{align}
Then there exist some positive number $\epsilon_0$ 
and an increasing function $g(t)$ with $g(0) = 0$ and $g(t)$ $=$ $O (t \log (1/t)^{d+1})$
$(t \searrow 0)$
such that if $z_0$ $\in$ $reg ( \epsilon_0 )$, that is
for some positive number $R_0$ and positive integer $\lambda_0$ possibly depending on $z_0$,
\begin{equation}
\frac 1{R_0^d} \lint_{P_{g(R_0)} (z_0)\cap\overline{Q(T)}}
\mathbf{e}_\lambda (z)\, dz\; < \;\epsilon_0\end{equation}
implies
\begin{equation}
\sup_{z\in P_{R_0}(z_0)\cap\overline{Q(T)}}\mathbf{e}_\lambda (z) 
\;\le\; C\Bigl[\frac 1{R_0^2}\,+\,\Vert u_0\Vert_{C^2(\partial\Omega)}\Bigr]
\label{INEQ:PR}\end{equation}
as long as any $\lambda$ is more than or equal to $\lambda_0$.\end{Thm}
\begin{Def}{}
In the sequel we respectively shorten $\,\mathbf{sing}$ and $\mathbf{reg}$
by $\mathbf{sing} ( \epsilon_0 )$ and $\mathbf{reg} ( \epsilon_0 )$.
\end{Def}
So we list a few results directly obtained from Lemma 3.6.
\begin{Lem}{}\label{LEM:Cont-Time}
Pick up any point $z_0$ $\in$ $\mathbf{reg}$ and fix it.
On the parabolic cylinder $P_{R_0/2}(z_0)$ $\cap$ $\overline{Q(T)}$, we obtain
\begin{equation}
|u_\lambda (t,x)\,-\,u_\lambda (s,x)|\;\le\;\frac C{R_0} |t-s|^{1/2}\label{INEQ:Cont-Time}
\end{equation}
for any points $t$ and $s$ in $[\max(0,t_0-(R_0/2)^2), t_0 + (R_0/2)^2]$
and $x$ $\in$ ${B}_{R_0/2}(x_0)$ $\cap$ $\overline{\Omega}$ with $z_0$ $=$ $(t_0,x_0)$.
\end{Lem}
\begin{Thm}{\rm{(Singular Set).}}\label{THM:Singular}
The set of $\mathbf{sing}$ is a closed set and 
\begin{equation}
\mathcal{H}^{(d)} (\mathbf{sing})\,=\,0
\label{EQ:Hausdorff-Estimate}
\end{equation}
holds with respect to the parabolic metric.
\end{Thm}
\begin{Thm}{\rm{(Strong Convergence)}}
\label{THM:Strong-Convergence-Gradient}
A sequence of the GLHF converges strongly to the WHHF in 
$L^2 ((0,T)\,;\,H^{1,2} (\Omega\cap K\,;\,\mathbb{R}^{D+1}))$
for any compact set $K$ $\subset$ $\mathbb{R}^d$. \end{Thm}%

%% file: main_bhf.tex
%
%#! platex bhf
%
%\setcounter{chapternumber}{4}\setcounter{equation}{0}
\renewcommand{\theequation}%
{\thechapternumber.\arabic{equation}}
\section{\enspace Proof of Main Theorems}
\par
By utilizing a few ingredients and properties on the WHHF and the GLHF,
this chapter proves Theorem \ref{THM:Main-1}
and Theorem \ref{THM:Main-3} in Chapter \ref{SEC:Intro}.
\subsection{Proof of Theorem \ref{THM:Main-1}}
We discuss a partial boundary regularity on the WHHF constructed above:
All proof proceeds as in the same way as the ones of chapter 4 in K.Horihata\cite{horihata}.
On account of Theorem \ref{THM:Singular} and Theorem \ref{THM:Strong-Convergence-Gradient}, 
we assert that
\allowdisplaybreaks\begin{align}
%1
& \mathbf{sing}\,=\,\bigcap_{R>0}\bigl\{z_0\in{\overline{Q(T)}}
\,;\,\frac 1{2R^d}\kern-2.5pt\lint_{P_R(z_0)\cap\overline{Q(T)}}\kern-2.5pt
|\nabla u|^2\,dz\,\ge\,\epsilon_0 \,\bigr\}
\label{EQ:Sing-2}\end{align}
is closed with $\mathcal{H}^d(\mathbf{sing})$
where a number $\epsilon_0$ is a positive constant appeared in Theorem \ref{THM:PR}. 
The smoothness for our WHHF on $\mathbf{reg}$ holds if \begin{math}u_0\,\in\,C^\infty(\partial\Omega).\end{math}
We readily see that the WHHF satisfies a monotonical inequality (i) and
a reverse Poincar\'{e} inequality (ii)
from Theorem \ref{THM:Mon} combined with Theorem \ref{THM:Strong-Convergence-Gradient}
and Corollary \ref{COR:RPI}.$\qed$
\vskip 6pt\subsection{Proof of Theorem \ref{THM:Main-3}}\vskip 6pt
\par
On account of Theorem \ref{THM:PR}, we have only to show that 
for any point $z_0$ $\in$ $\overline{Q(T)}$
\begin{equation}
\lint_{P_r (z_0)\cap\overline{Q(T)}} |\nabla u|^2\, dz\;<\;\frac{\epsilon_0^2  r^d}2
\end{equation}
holds with the positive number $\epsilon_0$ appeared in \eqref{EQ:Sing-2}.
\par By R.~Hamilton \cite{hamilton}, we find that there exists a positive number 
\begin{math}\tau_0\;=\;\tau_0(\Vert u_0\Vert_{C^{1,\alpha}})\;(0<\alpha_0<1)\end{math}
such that $u$ is of $C^{2,\alpha}$ $([0,\tau_0)\times\overline{\Omega})$.
If we use a stereo-graphic projection given 
by the mapping $v$ $=$ $(v_i)$ $(i=1,2,\ldots,D)$ of
\begin{equation}
\begin{array}{rl}
%1
u^i\,&=\;\displaystyle\frac{2v^i}{1+|v|^2} \quad (i\,=\,1,2,\ldots,D),\quad
u^D\,=\;\displaystyle\frac{1-|v|^2}{1+|v|^2}
\end{array}\end{equation}
with $|v|^2$ $=$ $\displaystyle\sum_{i=1}^{D}$ $(v^i)^2$,
our equation \eqref{EQ:GLHF} becomes for $v^i$ $(i\,=\,1,2,\ldots,D)$
\begin{equation}
\dfrac {\partial}{\partial t}\frac{v^i}{1+|v|^2}
\,=\,\triangle\frac{v^i}{1+|v|^2} 
\,+\,4\Bigl|\frac{\nabla v}{1+|v|^2}\Bigr|^2
\frac{v^i}{1+|v|^2}\quad\mathrm{in}\; Q(T).
\label{EQ:HHF-SG}\end{equation}
\par
By multiplying \eqref{EQ:HHF-SG} by $v^i$, we arrive at
\allowdisplaybreaks\begin{align}
&\left(\dfrac {\partial}{dt} \,-\,\triangle\right) (|v|^2-1)
\,+\,\frac{\vert\nabla |v|^2\vert^2}{(1+|v|^2)^2}
\;=\;\frac{|v|^2-1}{|v|^2+1}|\nabla{v}|^2.
\end{align}
\par Thus we obtain
\begin{equation}
\left(\dfrac {\partial}{\partial t}\,-\,\triangle\right)((|v|^2\,-\,1)^{(0)})^2
\;\le\;2\Vert\nabla v\Vert_{L^\infty}((|v|^2\,-\,1)^{(0)})^2
\label{INEQ:HHF-SG}
\notag\end{equation}
which shows \begin{math}|v|^2\,-\,1\;\le\;0\end{math} 
on \begin{math}[0,\tau_0)\times\overline{\Omega},\end{math}
where
\begin{equation*}
((|v|^2\,-\,1)^{(0)})\;=\;
\begin{cases}
|v|^2\,-\,1 &\quad\mathrm{if}\quad |v|^2\,\ge\,1,\\
0 &\quad\mathrm{if}\quad |v|^2\,<\,1.
\end{cases}\end{equation*}
So using \eqref{EQ:HHF-D} again, we arrive at
\begin{equation}
1\,-\,|v|^2\;\ge\;1\,-\,|v(0)|^2
\label{INEQ:Max-P}
\end{equation}
on \begin{math}[0,\tau_0)\times\overline{\Omega}\end{math}
\par We readily see that there exists a positive integer $k_0$ 
such that
\allowdisplaybreaks\begin{align}
& \underset{P_{d_0/2^{k_0-1}}(z_0)}{\esssup |v|^2}
\,-\,\underset{P_{d_0/2^{k_0}}(z_0)}{\esssup |v|^2}\;<\;\epsilon_0^2
\label{INEQ:F}\end{align} because of
\allowdisplaybreaks\begin{align}
%1
&\sum_{k=1}^\infty\bigl(\underset{P_{d_0/2^{k-1}}(z_0)}{\esssup|v|^2}
\,-\,\underset{P_{d_0/2^{k}} (z_0)}{\esssup|v|^2}\bigr)
\;=\;\underset{P_{d_0} (z_0)}{\esssup |v|^2}\,<\,\infty.\notag\end{align}
Recalling \eqref{EQ:HHF-D}, we deduce
\begin{equation}
\lint_{P_{d_0/2^k}(z_0)}|\nabla u|^2\,dz
\;<\;C(v_0)\epsilon_0^2 (d_0/2^{k})^d\;\le\;\epsilon_0 (d_0/2^{k})^d
\label{INEQ:Small-Energy-Density}\end{equation}
for any integer $k$ $\ge$ $k_0$ as long as a positive number $\epsilon_0$ is less than $1/C(v_0)$.
\par Moreover from \eqref{INEQ:Small-Energy-Density}, 
an application on an $\epsilon$-regularity theory by Y.Chen. Y.Li and F.H.Lin ~\cite{chen-li-lin} 
in which we exploit \eqref{THM:Existence-HF} verifies
\begin{equation}
\lint_{P_r (z_0)}\vert\nabla v \vert^2\,dz\;<\;Cr^{d+\alpha_0}\quad (0<r<d_0/2^{k_0+1})
\label{INEQ:E-Reg}
\end{equation}
where $\alpha_0$ is a positive number less than $1$,
and is independent of a point $z_0$ and a positive number $r$.
Thus we arrive at
\begin{equation}
\lint_{t_1-(2r)^2}^{t_1-r^2}\, dt
\lint_{B_{g(r)}(z_0)}\vert\nabla v\vert^2 \,dx\;<\;\epsilon_0 r^d,\notag\end{equation}
where the function $g$ was defined in Theorem \ref{THM:PR}.
\par Finally we prove the whole domain regularity: 
Assume that a point $(t_1,x_1)$ is a least-time singular point of $v$.
Namely \begin{math}|\nabla u|\;=\;\infty\end{math} at $(t_1,x_1)$
and \begin{math}|\nabla u|\;<\;\infty\end{math} 
on \begin{math}(0,t_1)\times\overline{\Omega}.\end{math}
Thus
\begin{equation}
\lint_{t_1-(2R)^2}^{t_1-R^2} \, dt \lint_{B_{g(R) R}}|\nabla u|^2 \, dx
\;\ge\;\epsilon_0 R^{\alpha_0}
\end{equation} for any positive $R$.
This is a contradiction.
\par Consequently from Theorem \ref{THM:PR}, we verifies a boundedness of $|\nabla v|$
which leads to our assertion.$\qed$
\begin{Rem}\label{REM:2ndRegularityCriterion}
\eqref{INEQ:Max-P} indicates that $u$ is smooth 
if $u(P_r)$ is compactly contained in $\mathbb{S}^{D-1}\cap\{y_d\,>\,0\}$
on any parabolic sylinder \begin{math}P_r\,\subset\!\subset\, Q(T).\end{math} 
Moreover since \begin{math}x/|x|:\mathbb{B}^3\,\to\,\mathbb{S}^2\end{math} 
is a unique weakly harmonic heat flow, 
the one-sided condition is sharp for the regularity to hold.
\end{Rem}

%% file: ref_bhf.tex
%
%#! platex bhf
%